\documentclass[12pt]{article}
\usepackage{amsmath}
\usepackage{amssymb}

\usepackage[cp1251]{inputenc}
\usepackage[russian]{babel}
\newcommand{\il}[2]{\int\limits_{#1}^{#2}}

\newcommand{\ph}{\phantom{a}}
\newcommand{\phh}{\phantom{aaa}}

\newcommand{\sist}[2]{\left\{
\begin{array}{l}
{#1}\\
\ph\\
{#2}
\end{array}
\right.}

\begin{document}
MSC 34C99

\vskip 10pt
{\bf \centerline {Global solvability criteria for  quaternionic}

\centerline{Riccati equations}

\vskip 20pt \centerline {G. A. Grigorian}}
\centerline{\it Institute  of Mathematics NAS of Armenia}
\centerline{\it E -mail: mathphys2@instmath.sci.am}
\vskip 20 pt

Abstract. Some   global  existence criteria for  quaternionic Riccati equations are \linebreak established. Two of them are  used to prove a completely non conjugation theorem for solutions of linear systems of ordinary differential equations.


\vskip 20pt
{\bf 1. Introduction}. Let  $a(t),\phantom{a} b(t),\phantom{a} c(t)$ and $d(t)$  be continuous quaternionic valued  functions on $[t_0;+\infty)$, i.e.:
$a(t)\equiv a_0(t) + i a_1(t) + j a_2(t) + k a_3(t), \phantom{a} b(t)\equiv b_0(t) + i b_1(t) + j b_2(t) + k b_3(t), \phantom{a} c(t)\equiv c_0(t) + i c_1(t) + j c_2(t) + k c_3(t), \phantom{a} d(t)\equiv d_0(t) + i d_1(t) + j d_2(t) + k d_3(t),$
where  $a_n(t), \phantom{a} b_n(t) \phantom{a}, c_n(t), \phantom{a} d_n(t) \phantom{a} (n=\overline{0,3})$ are real valued continuous functions on  $[t_0;+\infty)$, $i,\phantom{a} j, \phantom{a} k$ are the imaginary unities satisfying the conditions
$$
i^2 = j^2 = k^2 = ijk = -1, \phantom{a} ij = - ji = k. \eqno(1.1)
$$
Consider the quaternionic Riccati equation
$$
q' + q a(t) q + b(t) q + q c(t) + d(t) = 0, \phantom{aaa} t\ge t_0. \eqno (1.2)
$$
Here $q = q(t)$  is the unknown continuously differentiable quaternionic valued function.
Currently, there is a growing interest in quaternionic  differential equations, in particular, in Eq. (1.2) in connection with their various applications (see e.g., [1 -4]).
Criteria for the existence of periodic (and, therefore, global) solutions of Eq. (1.2) with periodic coefficients were obtained in [5, 6].
Explicit global existence criteria  for complex  solutions of Eq. (1.2) in the case of its complex coefficients were obtained in [7].

In this paper some  global  existence criteria for scalar quaternionic Riccati equations are obtained. Two of them are  used to prove a completely non conjugation theorem for solutions of linear systems of ordinary differential equations.

\hskip 10 pt

{\bf 2. Auxiliary propositions}. Substituting  $q= q_0 - i q_1 - j q_2 - k q_3$ in  (1.2),  where   $q_0$  is the real  and $-q_n (\overline{1,3})$   are the imaginary parts of $q$,  and separating the real and imaginary parts we come to the following nonlinear system
$$
\left\{\begin{array}{l}
q_0' + a_0(t) q_0^2 + \{b_0(t) + c_0(t) + 2 [a_1(t) q_1 + a_2(t) q_2 + a_3(t) q_3]\} q_0 - \\
\phantom{aaaaaaaaaaaaaaaaaaaaaaaaaaaaaaaaaaaaaaaaaaaaaa} - P(t, q_1, q_2, q_3) = 0;\\
q_1' + a_1(t) q_1^2 + \{b_0(t) + c_0(t) + 2 [a_0(t) q_0 + a_2(t) q_2 + a_3(t) q_3]\} q_1 -
 \\
\phantom{aaaaaaaaaaaaaaaaaaaaaaaaaaaaaaaaaaaaaaaaaaaaaa} -  Q(t, q_0, q_2, q_3) = 0;\\
q_2' + a_2(t) q_2^2 + \{b_0(t) + c_0(t) + 2 [a_0(t) q_0 + a_1(t) q_1 + a_3(t) q_3]\} q_2 - \\
\phantom{aaaaaaaaaaaaaaaaaaaaaaaaaaaaaaaaaaaaaaaaaaaaaa} - R(t, q_0, q_1, q_3) = 0;\\
q_3' + a_3(t) q_3^2 + \{b_0(t) + c_0(t) + 2 [a_0(t) q_0 + a_1(t) q_1 + a_2(t) q_2]\} q_3 - \\
\phantom{aaaaaaaaaaaaaaaaaaaaaaaaaaaaaaaaaaaaaaaaaaaaaa} - S(t, q_0, q_1, q_2) = 0;\\
\end{array}
\right. \eqno (2.1)
$$
where
$$
P(t,q_1,q_2,q_3) \equiv a_0(t)[q_1^2 + q_2^2 + q_3^2] - (b_1(t) + c_1(t)) q_1 - (b_2(t) + c_2(t)) q_2 - (b_3(t) + c_3(t)) q_3 - d_0(t);
$$
$$
Q(t,q_0,q_2,q_3) \equiv a_1(t)[q_0^2 + q_2^2 + q_3^2] + (b_1(t) + c_1(t)) q_0 +(b_3(t) - c_3(t)) q_2 - (b_2(t) - c_2(t)) q_3 + d_1(t);
$$
$$
R(t,q_0,q_1,q_3) \equiv a_2(t)[q_0^2 + q_1^2 + q_3^2] + (b_2(t) + c_2(t)) q_0 - (b_3(t) - c_3(t)) q_1 + (b_1(t) - c_1(t)) q_3 + d_2(t);
$$
$$
S(t,q_0,q_1,q_2) \equiv a_3(t)[q_0^2 + q_1^2 + q_2^2] + (b_3(t) + c_3(t)) q_0 +(b_2(t) - c_2(t)) q_1 - (b_1(t) - c_1(t)) q_2 + d_3(t);
$$
$t\ge t_0$.
Consider the square matrices
$$
E\equiv \left(\begin{array}{lccr}
1 &0 &0& 0\\
0 &1 &0 &0\\
0& 0& 1& 0\\
0& 0& 0& 1
\end{array}
\right),
\phantom{aaaaaaa}
I\equiv \left(\begin{array}{rccr}
0& 1& 0& 0 \\
-1& 0& 0& 0\\
0& 0& 0& 1 \\
0& 0& -1& 0
\end{array}
\right),
$$
$$
J\equiv \left(\begin{array}{rccr}
0& 0& 1& 0 \\
0& 0& 0& -1\\
-1& 0& 0& 0 \\
0& 1& 0& 0
\end{array}
\right),
\phantom{aaaaa}
K\equiv \left(\begin{array}{rccr}
0& 0& 0& -1 \\
0& 0 &-1 &0\\
0& 1& 0& 0 \\
1& 0& 0& 0
\end{array}
\right).
$$
It is not difficult to check that
$I^2 = J^2 = K^2 = IJK = -E, \phantom{a} IJ = - JI = K$. Then by (1.1) there is an one to one correspondence between the quaternions  $m\equiv m_0 + i m_1 + j m_2 + k m_3$ and the matrices of the form  $M\equiv m_0E + m_1 I + m_2 J + m_3 K$:
$$
m\equiv m_0 + i m_1 + j m_2 + k m_3 \leftrightarrow M\equiv \left(\begin{array}{rccr}
m_0& m_1& m_2& -m_3 \\
-m_1& m_0 &-m_3 & -m_2\\
-m_2& m_3& m_0& m_1 \\
m_3& m_2& -m_1& m_0
\end{array}
\right)
\eqno (2.2)
$$
The matrix $M$ corresponding  to the quaternion $m$ by the rule (2.2) we will call the symbol of the quaternion $m$ and will denote by  $\widehat{m}$.

Let  $A(t), \phantom{a} B(t), \phantom{a} C(t)$  and  $D(t)$  be the symbols of  $a(t),\phantom{a} b(t), \phantom{a} c(t)$   and  $d(t)$   respectively. Consider the matrix Riccati equation
$$
Y' + Y A(t) Y + B(t) Y + Y C(t) + D(t) = 0, \phantom{aaa} t\ge t_0. \eqno (2.3)
$$
By (2.2)  the solutions  $q(t)$  of Eq. (1.2), existing on some interval  $[t_1;t_2) \phantom{a} (t_0 \le t_1 < t_2 \le + \infty)$, are connected wit solutions  $Y(t)$ of Eq. (2.3) by equalities
$$
\widehat{q(t)} = Y(t), \phantom{aaa} t\in [t_1;t_2), \phantom{aaa}\widehat{q(t_1)} = Y(t_1). \eqno (2.4)
$$
Along with Eq. (2.3) consider the system of matrix equations
$$
\left\{
\begin{array}{l}
\Phi' = C(t) \Phi + A(t) \Psi;\\
\phantom{aaa}\\
\Psi' = - D(t) \Phi - B(t) \Psi, \phantom{aaa}t \ge t_0.
\end{array}
\right.
\eqno (2.5)
$$
Here $\Phi\equiv \Phi(t), \phantom{a} \Psi \equiv \Psi(t)$ are the unknown continuously differentiable matrix functions of dimension   $4\times 4$ on  $[t_0;+\infty)$.
Let $Y_0(t)$  be a solution of Eq. (2.3) on  $[t_1;t_2)$.  The substitution
$$
\Psi = Y_0(t) \Phi, \phantom{a} t\in [t_1;t_2), \eqno (2.6)
$$
in  (2.5) leads to the system

$$
\left\{
\begin{array}{l}
\Phi' = [A(t) Y_0(t) + C(t)] \Phi;\\
\phantom{a}\\

[Y'_0(t) + Y_0(t) A(t) Y_0(t) + B(t) Y_0(t) + Y_0(t) C(t) + D(t)]\Phi =0 \phantom{aaa} t\in [t_1;t_2).
\end{array}
\right.
$$

\noindent
Therefore  $(\Phi_0(t), Y_0(t) \Phi_0(t))$ is a solution of the system (2.5) on $[t_1;t_2)$,     where  $\Phi_0(t)$  is a solution to the following matrix equation
$$
\Phi' = [A(t) Y_0(t) + C(t)] \Phi, \phantom{aaa} t\in [t_1;t_2). \eqno (2.7)
$$
Let $Y(t) \phantom{a} (q(t))$  be a solution to Eq. (2.3)
(to Eq. (1.2)) on  $[t_1;t_2)$.

{\bf Definition 2.1}.  {\it The set $[t_1;t_2)$  is called the maximum existence interval for the solution $Y(t)$  of Eq. (2.3) (for the solution  $q(t)$  of Eq. (1.2)), if $Y(t) \phantom{a} (q(t))$  cannot be continued to the right from $t_2$.}

{\bf Lemma 2.1.} {\it Let $Y(t)$ be a solution of Eq. (2.3) on $[t_1; t_2) \phantom{a}  (t_0 \le t_1 < t_2 < + \infty)$.  Then $[t_1;t_2)$ is not the maximum existence interval for  $Y(t)$  provided the function
$$
f_0(t) \equiv \int\limits_{t_0}^t tr [A(\tau) Y(\tau)] d\tau, \phantom{aaa} t\in [t_1;t_2),
$$
is bounded from below on  $[t_1;t_2)$.
}

Proof. Let  $\Phi(t)$ be a solution to the matrix equation
$$
\Phi' = [A(t)Y(t) + C(t)]\Phi, \phantom{aaa} t\in [t_1; t_2), \phantom{aaa}\mbox{with}
$$
$$
\det \Phi(t_1) \ne 0. \eqno (2.8)
$$
By (2.6) and  (2.7), $(\Phi(t), Y(t)\Phi(t))$  is a solution to the system (2.5) on $[t_1;t_2)$  which can be continued on  $[t_0;+\infty)$ as a solution $(\Phi(t), \Psi(t))$   of the system (2.5). According to the Liouville's formula  (see [8], p. 46, Theorem 1.2)  we have:
$$
\det \Phi(t) = \det \Phi(t_1) \exp\biggl\{\int\limits_{t_0}^t tr [A(\tau) Y(\tau) + C(\tau)] d\tau\biggr\}, \phantom{aaa} t\in [t_1;t_2).
$$
From here from the conditions of lemma and from (2.8) it follows that   $\det \Phi(t) \ne 0, \linebreak t\in~ [t_1; t_3)$, for some $t_3 > t_2$.
Then by (2.6) and  (2.7) the matrix function \linebreak $\widetilde{Y}(t) \equiv \Psi(t) \Phi^{-1}(t), \phantom{a} t\in [t_1;t_3),$ is a solution to Eq. (2.3) on  $[t_1;t_3)$. Obviously  $\widetilde{Y}(t)$  coincides with $Y(t)$  on  $[t_1;t_2)$.   Therefore  $[t_1;t_2)$  is not the maximum existence interval for $Y(t)$. The lemma is proved.

Let   $f(t), \phantom{a} g(t), \phantom{a} h(t), \phantom{a} f_1(t), \phantom{a} g_1(t),  \phantom{a}h_1(t)$  be real valued continuous functions on \linebreak $[t_0;+\infty)$.
Consider the Riccati equations
$$
y' + f(t) y^2 + g(t) y + h(t) = 0, \phantom{aaa} t\ge t_0; \eqno (2.9)
$$
$$
y' + f_1(t) y^2 + g_1(t) y + h_1(t) = 0, \phantom{aaa} t\ge t_0. \eqno (2.10)
$$
and the differential inequalities
$$
y' + f(t) y^2 + g(t) y + h(t) \ge 0, \phantom{aaa} t\ge t_0; \eqno (2.11)
$$
$$
y' + f_1(t) y^2 + g_1(t) y + h_1(t) \ge 0, \phantom{aaa} t\ge t_0. \eqno (2.12)
$$

{\bf Remark 2.1} {\it  For   $f(t) \ge 0, \phantom{a} t\ge t_0,$ every solution of the linear equation  \linebreak $y' + g(t) y + h(t) = 0$  on  $[t_0;\tau_0) \phantom{a} (t_0 < \tau_0 \le +\infty)$ is a solution of the inequality  (2.11) on $[t_0;\tau_0)$}.

{\bf Remark 2.2} {\it Every solution of Eq. (2.10) on $[t_0;\tau_0) \phantom{a} (t_0 < \tau_0 \le +\infty)$ is also a solution of the inequality (2.12) on  $[t_0;\tau_0)$}.

{\bf Theorem 2.1}. {\it  Let Eq. (2.10) has a real solution  $y_1(t)$  on $[t_0;\tau_0) (\tau_0 \le +\infty)$,  and let the following conditions be satisfied: $f(t) \ge 0$,
$$
\int\limits_{t_0}^t\exp\biggl\{\int\limits_{t_0}^\tau \bigl[ f(s)(\eta_0(s) + \eta_1(s)) + g(s)\bigr] d s\biggr\}\times \phantom{aaaaaaaaaaaaaaaaaaaaaaaaaaaaaaaaaaaaaaa}
$$
$$
\phantom{aaa}\times\Bigl[(f_1(\tau) - f(\tau)) y_1^2(\tau) + (g_1(\tau) - g(\tau)) y_1(\tau) + h_1(\tau) - h(\tau)\Bigr] d\tau \ge 0, \phantom{aaa} t\in [t_0;\tau_0).
$$
where  $\eta_0(t)$   and  $\eta_1(t)$ are solutions of the inequalities (2.11) and  (2.12) on  $[t_0;\tau_0)$ such that $\eta_j(t_0) \ge y_1(t_0), \phantom{a} j=0,1$.
Then for every $\gamma_0 \ge y_1(t_0)$ Eq. (2.9) has a solution $y_0(t)$   on $[t_0;\tau_0),$ satisfying the initial conditions  $y_0(t_0) = \gamma_0$,  and $y_0(t) \ge y_1(t), \phantom{aa} t \in [t_0;\tau_0).$}

The proof of this theorem is presented in [9].

Let $t_0 < t_1 < ... $ be a finite or infinite sequence such that $t_m \in [t_0;\tau_0] \phantom{a} (t_0 < \tau_0 \le + \infty)$. We assume that if $\{t_m\}$ is finite then $\max \{t_m\} = \tau_0$ otherwise $\lim\limits_{m \to +\infty} t_m = \tau_0$. Denote: $I_{g,h}(\xi,t) \equiv \int\limits_\xi^t\exp\biggl\{- \int\limits_\tau^tg(s) d s\biggr\} h(\tau) d \tau, \phantom{a} t \ge \xi \ge t_0.$

{\bf Theorem 2.2.} {\it Let $f(t) \ge 0, \phantom{a} t\in [t_0;\tau_0),$ and
$$
\int\limits_{t_k}^t \exp\biggl\{\int\limits_{t_k}^\tau\Bigl[g(s) - f(s) I_{g,h}(t_k, s)\Bigr] d s\biggr\} h(\tau) d \tau \le 0, \phantom{a} t\in [t_k;t_{k+1}), \phantom{a} k =1, 2, ... .
$$
Then  for every $\gamma_0 \ge 0$ Eq. (2.9) has a solution $y_0(t)$ on $[t_0;\tau_0)$ satisfying the initial condition $y_0(t_0) = \gamma_0$ and $y_0(t) \ge 0, \phantom{a} t\in [t_0;\tau_0)$.}

See the proof in [10]

{\bf Theorem 2.3}. {\it Let $\alpha(t)$ and $\beta(t)$ be continuously differentiable on $[t_0; \tau_0)$ functions and
$
\alpha(t) > 0, \phantom{aa} \beta(t) > 0, \phantom{a} t\in [t_0;\tau_0);
$

\noindent
A) $0 \le f(t) \le \alpha(t), \phantom{a} h(t) \le \beta(t), \phantom{a} t\in [t_0;\tau_0);$

\noindent
B) $g(t) \ge \frac{1}{2} \biggl[\frac{\alpha'(t)}{\alpha(t)} - \frac{\beta'(t)}{\beta(t)}\biggr] + 2\sqrt{\alpha(t) \beta(t)}, \phantom{aaa} t\in [t_0; \tau_0).$

\noindent
Then  for every $\gamma_0 \ge - \sqrt{\frac{\beta(t_0)}{\alpha(t_0)}}$  Eq. (2.9) has a solution $y_0(t)$ on $[t_0;\tau_0)$ with $y_0(t_0) = \gamma_0$ and
$$
y_0(t) \ge - \sqrt{\frac{\beta(t)}{\alpha(t)}}, \phantom{aaa} t\in [t_0; \tau_0).
$$
}

See the proof in ]11, Theorem 8].

{\bf Theorem 2.4.} {\it Let $\alpha(t)$ and $\beta(t)$ be the same as in Theorem 2.3. If assumption A) of Theorem 2.3 and the inequality

\noindent
D)  $g(t) \le \frac{1}{2} \biggl[\frac{\alpha'(t)}{\alpha(t)} - \frac{\beta'(t)}{\beta(t)}\biggr] - 2\sqrt{\alpha(t) \beta(t)}, \phantom{aaa} t\in [t_0; \tau_0),$

\noindent
are valid, then  for every $\gamma_0 \ge \sqrt{\frac{\beta(t_0)}{\alpha(t_0)}}$ Eq. (2.9) has a solution $y_0(t)$ on $[t_0;\tau_0)$ with $y_0(t_0) =~ \gamma_0$ and
$$
y_0(t) \ge  \sqrt{\frac{\beta(t)}{\alpha(t)}}, \phantom{aaa} t\in [t_0; \tau_0).
$$
}

See the proof in ]11, Theorem 7].

{\bf Theorem 2.5}. {\it Let $\alpha_m(t)$ and $\beta_m(t), \ph m=1,2$, be continuously differentiable functions
 on $[t_0;\tau_0)$, and let  $(-1)^m \alpha_m(t) > 0, \ph (-1)^m \beta_m(t) > 0, \ph t\in[t_0;\tau_0), \ph m=1,2,$.
 If:

\noindent
E) $\alpha_1(t) \le f(t) \le \alpha_2(t),\ph \beta_1(t) \le h(t) \le \beta_2(t), \ph t\in [t_0;\tau_0);$

\noindent
F) $g(t) \ge  \frac{1}{2}\biggl(\frac{\alpha'_m(t)}{\alpha_m(t)} - \frac{\beta'_m(t)}{\beta_m(t)}\biggr) +
 2(-1)^m \sqrt{\alpha_m(t)\beta_m(t)}, \ph t\in [t_0;\tau_0), \ph m=1,2,$

\noindent
then for any $y_{(0)} \in \biggl[-\sqrt{\frac{\beta_2(t_0)}{\alpha_2(t_0)}};
 \sqrt{\frac{\beta_1(t_0)}{\alpha_1(t_0)}}\biggr]$  Eq. (2.9) has a solution $y_0(t)$ on $[t_0;\tau_0)$
  satisfying the initial condition  $y_0(t_0) = y_{(0)}$, and
$$
-\sqrt{\frac{\beta_2(t)}{\alpha_2(t)}} \le y_0(t) \le \sqrt{\frac{\beta_1(t)}{\alpha_1(t)}}, \phh t\in [t_0;\tau_0).
$$
}

See the proof  in [10, Theorem 4.2] .

Let $p, \phantom{a} q, \phantom{a} r,  \phantom{a} s,  \phantom{a} l$ be real numbers and let $\varepsilon > 0$.

{\bf Definition 2.2.} {\it  The ordered fiver $(p, q, r, s, l)$  is called  $\varepsilon$ - semi definite positive if:

\noindent
1) \phantom{a} $p> 0, \phantom{aa} l > 0$;

\noindent
2) $\max\{q, r, s\} \ge \sqrt{l + \varepsilon}$  or

\noindent
$0 \le \min\{q, r, s\} \le \max\{q, r, s\} \le \sqrt{l + \varepsilon}$ and

\noindent
$q^2 + r^2 + s^2 \ge l+s$.}

{\bf Remark 2.3.}  {\it From the geometrical point of view the relations 1)  and  2) mean that the ball of radius  $\sqrt{l+ \varepsilon}$ with its center in the point $(q,r,s)$  may be located in any such position in the space of coordinates $x, \phantom{a} y, \phantom{a} z$, that its intersection with  the octant  $x >~ 0, \phantom{a} y >~ 0, \phantom{a} z > 0$ is empty.}

Consider the quadratic form
$$
W(x,y,z)\equiv p\Bigl[\Bigl(x + \frac{q}{2p}\Bigr)^2 + \Bigl(y + \frac{r}{2p}\Bigr)^2  + \Bigl(z + \frac{s}{2p}\Bigr)^2\Bigr]  - \frac{l}{4p}, \phantom{aaa} x, y, z \in (-\infty; + \infty).
$$

 {\bf Lemma 2.2}. {\it  If for some $\varepsilon > 0$   the ordered fiver   $(p, q, r, s, l)$ is
  $\varepsilon$ - semi definite positive then for every  $x\ge 0, \phantom{a} y\ge 0, \phantom{a} z\ge 0$ the
 inequality
$$
W(x,y,z) \ge \varepsilon/4p
$$                                                                                                                      is satisfied.}

 Proof. For every $x\ge 0, \phantom{a} y\ge 0, \phantom{a} z\ge 0 \phantom{a}$ we have: if  $\max\{q, r, s\} \ge
 \sqrt{l + \varepsilon}$,   then  $W(x, y, z) \ge p\frac{l + \varepsilon}{4 p^2} - \frac{l}{4 p}= \frac{\varepsilon}{4
 p}$, and if  $ 0 \le \min\{q, r, s\} \le \max\{q, r, s\} \le \sqrt{l + \varepsilon}$, then since $q \ge 0, \phantom{a}
 r\ge 0, \phantom{a} s\ge 0$,we will get: $W(x, y, z) \ge p\bigl(\frac{q^2}{4 p^2} + \frac{r^2}{4 p^2} + \frac{s^2}{4
 p^2}\bigr)
  - \frac{l}{4 p} \ge \frac{l + \varepsilon}{4 p} - \frac{l}{4 p} = \frac{\varepsilon}{4 p}$. The lemma is proved.

 \vskip 10pt

{\bf 3. Global solvability criteria}. In this section we  study the global solvability conditions of Eq. (1.2) in the case when
 $a_n(t) \ge 0, \phantom{a} t\ge t_0, \phantom{a} n = \overline{0,3}$.
The cases when $(-1)^{m_n} a_n(t) \ge~ 0, \linebreak t\ge t_0, \phantom{a} m_n = 0,1, \phantom{a} n=\overline{0,3}, \phantom{a} m_0 + m_1 + m_2 + m_3 > 0$
are reducible to the studying one by the simple transformations  $q \rightarrow -q, \phantom{a} q \rightarrow  \overline{q}, \phantom{a} q \rightarrow iq, \phantom{a} q \rightarrow jq, \ph q \rightarrow kq$ and their combinations
 in  (1.2).   Denote:

\noindent
$p_{0,m}(t) \equiv b_m(t) + c_m(t), \phantom{a} m = \overline{1,3}, \phantom{a}  p_{1,1}(t) \equiv b_1(t) + c_1(t), \phantom{a} p_{1,2}(t) \equiv b_2(t) - c_2(t), \phantom{a} p_{1,3}(t) \equiv b_3(t) - c_3(t), \phantom{a} p_{2,1}(t) \equiv b_1(t) - c_1(t), \phantom{a} p_{2,2}(t) \equiv b_2(t) + c_2(t), \phantom{a} p_{2,3}(t) \equiv b_3(t) - c_3(t), \phantom{a} p_{3,m}(t) \equiv b_m(t) - c_m(t), \phantom{a} m = \overline{1,3}, \phantom{a} t\ge t_0. $
$$
D_0(t)\equiv\left\{\begin{array}{l}
\sum\limits_{m=1}^3 p_{0,m}^2(t) + 4 a_0(t) d_0(t), \phantom{a} if \phantom{a} a_0(t) \ne 0;\\
4 d_0(t) \phantom{a} if \phantom{a} a_0(t) = 0,
\end{array}\right. \phantom{aaaaaaaaaaaaaaaaaaaaaaaaaaaaaaaaaaaaa}
$$
$$
\phantom{aaaaaaa} D_n(t)\equiv\left\{\begin{array}{l}
\sum\limits_{m=1}^3 p_{n,m}^2(t) - 4 a_n(t) d_n(t), \phantom{a} if \phantom{a} a_n(t) \ne 0;\\
-4 d_n(t) \phantom{a} if \phantom{a} a_n(t) = 0,
\end{array}\right. \phantom{aaa} n = \overline{1,3}, \phantom{aaa} t\ge t_0.
$$
Let $\mathfrak{S}$ be a nonempty  subset of the set $\{0, 1, 2, 3\}$ and let $\mathfrak{O}$ be its complement, i.e., $\mathfrak{O} = \{0, 1, 2, 3\} \backslash \mathfrak{S}$

{\bf Theorem 3.1}. {\it Assume  $a_n(t) \ge 0, \phantom{a}  n \in \mathfrak{S}$ and if $a_n(t) = 0$  then $p_{n,m}(t) = 0, \phantom{a} m =~ \overline{1,3}, \linebreak   n\in \mathfrak{S}; \phantom{a} a_n(t) \equiv 0, \phantom{a}  n \in  \mathfrak{O}, \phantom{a} D_n(t) \le 0, \phantom{a} n \in \mathfrak{S}, \phantom{a} t \ge t_0.$

Then for every  $\gamma_n \ge 0, \phantom{a} n \in \mathfrak{S}, \phantom{a} \gamma_n \in (- \infty; + \infty), \phantom{a} n \in \mathfrak{O}$,  Eq. (1.2) has a solution  $q(t) \equiv q_0(t) - i q_1(t) - j q_2(t) - k q_3(t)$ on $[t_0; +\infty)$  with  $q_n(t_0) = \gamma_n, \phantom{a} n = \overline{0,3}$  and
$$
q_n(t) \ge 0, \phantom{aaa} n \in \mathfrak{S}, \phantom{aaa}   t\ge t_0. \eqno (3.1)
$$
Moreover if for some $n \in \mathfrak{S}$  $\gamma_n > 0$,  then also  $q_n(t) > 0$.}

Proof. Let  $[t_0;T)$ be the maximum existence interval for the solution  $q(t) \equiv q_0(t) - i q_1(t) - j q_2(t) - k q_3(t)$  of Eq. (1.2) satisfying the initial conditions $q_n(t_0) = \gamma_n, \phantom{a} n =~\overline{0,3}$ (existence of $[t_0;T)$   follows from the theory of normal systems of ordinary differential equations and from (2.1)). Show that
$$
q_n(t) \ge 0, \phantom{aaa} t\in [t_0; T), \phantom{aaa} n \in   \mathfrak{S}. \eqno (3.2)
$$
Let us prove the theorem in the case when $0 \in  \mathfrak{S}$. The proof of the theorem for other nonempty $\mathfrak{S}$ can be proved by analogy.
Consider the Riccati equations
$$
x' + a_0(t) x^2 + \{b_0(t) + c_0(t) + 2 [a_1(t) q_1(t) + a_2(t) q_2(t) + a_3(t) q_3(t)]\} x -\phantom{aaaaaaaaaaaaaaaaaaaaaa}
$$
$$
\phantom{aaaaaaaaaaaaaaaaaaaaaaaaaaaaaaaaa}- P(t,q_1(t), q_2(t), q_3(t)) = 0, \phantom{aa} t\in [t_0;T), \eqno (3.3)
$$
$$
x' + a_0(t) x^2 + \{b_0(t) + c_0(t) + 2 [a_1(t) q_1(t) + a_2(t) q_2(t) + a_3(t) q_3(t)]\} x = 0, \phantom{aa} t\in [t_0;T). \eqno (3.4)
$$
From the conditions of the theorem it follows that  $P(t, q_1(t), q_2(t), q_3(t)) \ge 0, \phantom{a} t\in [t_0;T)$. Then using Theorem 2.1
to the equations (3.3) and  (3.4) we conclude that the solution  $x(t)$
of Eq. (3.3) with  $x(t_0) = \gamma_0 \ge 0$ exists on $[t_0;T)$ and is non negative  (since  $x_1(t) \equiv 0$ is a solution to Eq. (3.4) on  $[t_0;T)$). Obviously $q_0(t)$ is a solution of Eq. (3.3). Hence $q_0(t) = x(t) \ge 0, \ph t\in [t_0; T)$. By analogy can be proved the remaining inequalities  (3.2). By (2.4) $Y(t) \equiv \widehat{q(t)}, \phantom{a} t\in [t_0;T),$  is a solution of Eq. (2.3) on   $[t_0;T)$. Then it is not difficult to verify that $tr [A(t)Y(t)] = \sum\limits_{n=0}^3 a_n(t) q_n(t) =  \sum\limits_{n\in \mathfrak{S}} a_n(t) q_n(t), \phantom{a} t\in [t_0;T).$  From here and from (3.2) we have:
$$
tr [A(t)Y(t)] \ge 0, \phantom{aaa} t\in [t_0;T). \eqno (3.5)
$$
Show that
$$
T = + \infty. \eqno (3.6)
$$
Suppose $T < + \infty$.  Then by virtue of Lemma 2.1 from (3.5) it follows that $[t_0;T)$  is not the maximum existence interval for $Y(t)$. Therefore  $[t_0;T)$  is not the maximum existence interval for $q(t)$.  The obtained contradiction proves (3.6).
From (3.6) and  (3.2) it follows (3.1). Assume $\gamma_0 > 0$.  By already  proven the solution $\widetilde{x}(t)$ of Eq. (3.3)  with  $\widetilde{x}(t_0) = 0$  exists on $[t_0; +\infty)$  and is nonnegative. Then by virtue of Theorem 2.1 the solution  $x(t)$  of Eq. (3.3) with   $x(t_0) = \gamma_0 > 0$ exists on  $[t_0; + \infty)$ and  $x(t) \ne \widetilde{x}(t), \phantom{a} t\ge t_0$. Therefore  $x(t) > 0, \phantom{a} t\ge t_0$. Obviously  $x(t) \equiv q_0(t), \phantom{a} t\ge t_0$. Hence  $q_0(t) > 0, \phantom{a} t\ge t_0$.    By analogy it can be shown that if  $\gamma_n > 0$  for some other  $n \in \mathfrak{S}$,  then also  $q_n(t) > 0, \phantom{a} t\ge t_0.$ The theorem is proved.

{\bf Remark 3.1.} {\it Theorem 3.1 is a generalization of Theorem 3.1 of work [7]}.

Set:
$
\mathcal{L}_0(t)\equiv (a_0(t), -b_1(t) - c_1(t), -b_2(t) - c_2(t), -b_3(t) - c_3(t),\phantom{a} D_0(t));
$

$
\phantom{aaaaaa} \mathcal{L}_1(t)\equiv (a_1(t),\phantom{a} b_1(t) + c_1(t), -b_2(t) + c_2(t),\phantom{a} b_3(t) - c_3(t),\phantom{a} D_1(t))
$

$
\phantom{aaaaaaaa}\mathcal{L}_2(t)\equiv(a_2(t),\phantom{a} b_1(t) - c_1(t),\phantom{a} b_2(t) + c_2(t),\phantom{a} b_3(t) - c_3(t),\phantom{a} D_2(t))
$

$
\phantom{aaaaaaaaaa}\mathcal{L}_3(t)\equiv (a_3(t), -b_1(t) + c_1(t),\phantom{a} b_2(t) - c_2(t),\phantom{a} b_3(t) + c_3(t),\phantom{a} D_3(t))
$

{\bf Theorem 3.2.} {\it  Let
 for some  $\varepsilon > 0$
and for every  $t\ge t_0$   the ordered fivers  $\mathcal{L}_n(t), \linebreak n=\overline{0,3}$
be  $\varepsilon$ -  semi definite  positive. Then for every  $\gamma_n > 0, \phantom{a} n =~\overline{0,3},$ Eq. (1.2) has a solution  $q(t) \equiv q_0(t) - i q_1(t) - j q_2(t) - k q_3(t)$  on $[t_0; +\infty)$ with  $q_n(t_0) = \gamma_n,  \phantom{a} n = \overline{0,3}$,  and
$$
q_n(t) > 0, \phantom{aa} t\ge t_0, \phantom{aaa} n=\overline{0,3}. \eqno (3.7)
$$
}

Proof. Let $[t_0;T)$ be the maximum existence interval for the solution  $q(t) \equiv q_0(t) - i q_1(t) - j q_2(t) - k q_3(t)$     of Eq. (1.2) satisfying the initial conditions  $q_n(t_0) = \gamma_n  \phantom{a} n = \overline{0,3}$.  Show that
$$
q_n(t) \ge 0, \phantom{aa} t\in [t_0;T) \phantom{aaa} n = \overline{0,3}. \eqno (3.8)
$$
Set:  $T_1 \equiv \sup\{t\in [t_0;T) : q_n(t) \ge 0, \phantom{a} t\in [t_0;T) \phantom{a} n = \overline{0,3}\}.$ Suppose (3.8) is not true. Then (obviously  $T_1 > t_0)$
$$
T_1 < T. \eqno (3.9)
$$
On the other hand from the conditions of the theorem it follows that

$
P(t,\phantom{a} q_1(t),\phantom{a} q_2(t),\phantom{a} q_3(t)) \ge \frac{\varepsilon}{4 a_0(t)},\phantom{aaa} Q(t,\phantom{a} q_0(t),\phantom{a} q_2(t),\phantom{a} q_3(t)) \ge \frac{\varepsilon}{4 a_1(t)},\phantom{aaa}
$

\phantom{aaa}
$
R(t,\phantom{a} q_0(t),\phantom{a} q_1(t),\phantom{a} q_3(t)) \ge \frac{\varepsilon}{4 a_2(t)},\phantom{aaa} S(t,\phantom{a} q_0(t),\phantom{a} q_1(t),\phantom{a} q_2(t)) \ge \frac{\varepsilon}{4 a_3(t)}, \phantom{aa} t\in [t_0;T_1).
$

\noindent
By the continuity property of the functions  $P, Q, R, S, q_0, q_1, q_2$ and $q_3$   it follows that for some  $T_2 > T_1 \phantom{a} (T_2 < T)$ the following inequalities are fulfilled:
$$
\left\{
\begin{array}{l}
P(t,\phantom{a} q_1(t),\phantom{a} q_2(t),\phantom{a} q_3(t)) \ge 0; \phantom{aaa} Q(t,\phantom{a} q_0(t),\phantom{a} q_2(t),\phantom{a} q_3(t)) \ge 0,\\
\phantom{a}\\
R(t,\phantom{a} q_0(t),\phantom{a} q_1(t),\phantom{a} q_3(t)) \ge 0; \phantom{aaa} S(t,\phantom{a} q_0(t),\phantom{a} q_1(t),\phantom{a} q_2(t)) \ge 0,
\end{array}
\right.
\eqno (3.10)
$$
For all $t\in [t_0;T_2)$. Consider on  $[t_0;T_2)$  the Riccati equations
$$
x' + a_0(t) x^2 + \{b_0(t) + c_0(t) + 2[a_1(t) q_1(t) + a_2(t) q_2(t) + a_3(t) q_3(t)]\} x - \phantom{aaaaaaaaaaaaaaaaaa}
$$
$$
\phantom{aaaaaaaaaaaaaaaaaaaaaaaaaaaaaaaaaaaa} - P(t, q_1(t), q_2(t), q_3(t)) = 0; \eqno (3.11)
$$
$$
x' + a_1(t) x^2 + \{b_0(t) + c_0(t) + 2[a_0(t) q_0(t) + a_2(t) q_2(t) + a_3(t) q_3(t)]\} x - \phantom{aaaaaaaaaaaaaaaaaa}
$$
$$
\phantom{aaaaaaaaaaaaaaaaaaaaaaaaaaaaaaaaaaaa} - Q(t, q_0(t), q_2(t), q_3(t)) = 0; \eqno (3.12)
$$
$$
x' + a_2(t) x^2 + \{b_0(t) + c_0(t) + 2[a_0(t) q_0(t) + a_1(t) q_1(t) + a_3(t) q_3(t)]\} x - \phantom{aaaaaaaaaaaaaaaaaa}
$$
$$
\phantom{aaaaaaaaaaaaaaaaaaaaaaaaaaaaaaaaaaaa} - R(t, q_0(t), q_1(t), q_3(t)) = 0; \eqno (3.13)
$$
$$
x' + a_3(t) x^2 + \{b_0(t) + c_0(t) + 2[a_0(t) q_0(t) + a_1(t) q_1(t) + a_2(t) q_2(t)]\} x - \phantom{aaaaaaaaaaaaaaaaaa}
$$
$$
\phantom{aaaaaaaaaaaaaaaaaaaaaaaaaaaaaaaaaaaa} - S(t, q_0(t), q_1(t), q_2(t)) = 0. \eqno (3.14)
$$

Let  $x_0(t), \phantom{a}  x_1(t), \phantom{a}  x_2(t)$  and $x_3(t)$   be the solutions of the equations (3.11),  (3.12), (3.13)  and  (3.14) respectively with  $x_n(t_0) = 0, \phantom{a} n = \overline{0,3}$. By virtue of Theorem 2.1 from (3.10) it follows that $x_n(t), \phantom{a} n = \overline{0,3}$, exist on  $[t_0;T_2)$ and are non negative. Then since \linebreak $q_0(t), \phantom{a}  q_1(t), \phantom{a}  q_2(t)$  and $q_3(t)$ are solutions of the equations (3.11),  (3.12), (3.13)  and (3.14) on $[t_0;T_2)$    and  $q_n(t_0) > x_n(t_0), \phantom{a} n = \overline{0,3}$,
the last functions (i.e. $q_n(t), \ph n=\overline{0,3}$)  are also non negative on  $[t_0;T_2)$,
which contradicts (3.9). The obtained contradiction proves (3.8). By virtue of Lemma 2.2 from  (3.8) and from the conditions of the theorem it follows that on $[t_0;T)$  the inequalities (3.10) are fulfilled.  Hence the solutions $x_n(t_0) \phantom{a} (n = \overline{0,3})$ exist on $[t_0;T)$ and are non negative. Obviously  $q_0(t), \phantom{a}  q_1(t), \phantom{a}  q_2(t)$  and $q_3(t)$  are solutions of the equations (3.11), (3.12), (3.13)  and (3.14) respectively on  $[t_0;T)$ and  $q_n(t_0) > x_n(t_0), \phantom{a} n = \overline{0,3}$. Therefore  $q_n(t) > 0, \phantom{a} t\in [t_0;T), \phantom{a} n = \overline{0,3}$.     Further, the proof of the theorem is  carried out similar to the proof of Theorem 3 1. The theorem is proved.

{\bf Theorem 3.3.} {\it Let $a_0(t) \ge 0, \phantom{a} a_n(t)\equiv 0, \phantom{a} n= \overline{1,3}, \phantom{a} t\ge t_0$, and
$$
\int\limits_{t_m}^t \exp\biggl\{\int\limits_{t_m}^t\bigl[b_0(s) + c_0(s)  - I_{b_0+ c_0, D_0}(t_m,s)\bigr] d s\biggr\} D_0(\tau) d \tau \le 0, \phantom{a} t\in [t_m; t_{m+1}), \phantom{a} m =0, 1, ....
$$
Then for every $\gamma_0 \ge 0, \phantom{a} \gamma_n \in (- \infty; + \infty), \phantom{a} n = \overline{1,3},$ Eq. (1.2) has a solution $q(t) \equiv q_0(t) - i q_1(t) - j q_2(t) - k q_3(t)$ with $q_n(t_n) = \gamma_n, \phantom{a} n = \overline{0,3},$ on $[t_0; +\infty)$ and
$$
q_0(t) \ge 0, \phantom{aaa} t\ge  t_0.  \eqno (3.15)
$$
}

Proof. Let $q(t) \equiv q_0(t) - i q_1(t) - j q_2(t) - k q_3(t)$ be the solution of Eq. (1.2) with $q_n(t_0) = \gamma_n, \phantom{a} n = \overline{0,3},$ and let $[t_0; T)$ be the maximum existence interval for $q(t)$. Show that
$$
T  = + \infty. \eqno (3.16)
$$
Consider the Riccati equation
$$
y' + a_0(t) y^2 + [b_0(t) + c_0(t)] y + D_0(t) = 0, \phantom{aaa} t\ge t_0. \eqno (3.17)
$$
By Theorem 2.2 from the conditions of the theorem it follows that for every $\gamma_0 \ge 0$ this equation has a solution $y_0(t)$ on $[t_0; + \infty)$ and $y_0(t) \ge 0, \phantom{a} t\ge t_0$. Then using Theorem~ 2.1 to  Eq. (3.11)
 and Eq. (3.17) and taking into account the fact that $q_0(t)$ is a solution to Eq. (3.11) we conclude that
$$
q_0(t) \ge y_0(t) \ge 0, \phantom{aaa} t\ge t_0. \eqno (3.18)
$$
Suppose $T < +\infty$. Then from (3.18) it follows that
$$
tr[A(t) Y(t)] = \int\limits_{t_0}^t a_0(s) q_0(s) d s \ge 0, \phantom{aaa} t\in [t_0;T).
$$
By virtue of Lemma 2.1 from here it follows that $[t_0;T)$ is not the maximum existence interval for $q(t)$ which contradicts our supposition. The obtained contradiction proves (3.16). From (3.16) and (3.18)
 it follows (3.15). The theorem is proved.

{\bf Remark 3.2.} {\it Unlike of the conditions of Theorem 3.1 and Theorem 3.2 the conditions of  Theorem~ 3.3
allow  $D_0(t)$ to change sign in every $[t_m;t_{m+1}),\ph m = 0, 1, ...$.}

By use of Theorem 2.3 and Theorem 2.4 analogically can be proved the following two theorems respectively

{\bf Theorem 3.4}. {\it Let $\alpha(t)$ and $\beta(t)$ be continuously differentiable on $[t_0; +\infty)$ functions and $\alpha(t) > 0, \phantom{a} \beta(t) > 0, \phantom{a} t\ge t_0,$

\noindent
$A_1) \phantom{a} 0\le a_0(t) \le \alpha(t), \phantom{a} D_0(t) \le \beta(t), \phantom{a} a_n(t) \equiv 0, \phantom{a} n = \overline{1,3}, \phantom{a} t\ge t_0 ;$

\noindent
$B_1) \phantom{a} b_0(t) + c_0(t) \ge \frac{1}{2}\biggl[\frac{\alpha'(t)}{\alpha(t)} - \frac{\beta'(t)}{\beta(t)}\biggr] + \sqrt{\alpha(t) \beta(t)}, \phantom{a} t \ge t_0.$

\noindent
Then for every $\gamma_0 \ge - \sqrt{\frac{\beta(t_0)}{\alpha(t_0)}}, \phantom{a} \gamma_n \in (-\infty; + \infty), \phantom{a} n = \overline{1,3},$ Eq. (1.2) has a solution $q(t)\equiv q_0(t) - i q_1(t) - j q_2(t) - k q_3(t)$ on $[t_0;+ \infty)$ with $q_n(t_0) = \gamma_n, \phantom{a} n = \overline{0,3}$, and
$$
q_0(t) \ge - \sqrt{\frac{\beta(t)}{\alpha(t)}}, \phantom{aaa} t\ge t_0.
$$
\phantom{aaaaaaaaaaaaaaaaaaaaaaaaaaaaaaaaaaaaaaaaaaaaaaaaaaaaaaaaaaaaaaaaaaaaaaaaa} $\Box$}

{\bf Theorem 3.5}. {\it Let $\alpha(t)$ and $\beta(t)$ be the same as in Theorem 3.4. If assumption $A_1)$ of Theorem 3.4 and the inequality

\noindent
$C_1) \phantom{a} b_0(t) + c_0(t) \le  \frac{1}{2}\biggl[\frac{\alpha'(t)}{\alpha(t)} - \frac{\beta'(t)}{\beta(t)}\biggr] - \sqrt{\alpha(t) \beta(t)}, \phantom{a} t\ge t_0,$

\noindent
are valid. Then for every  $\gamma_0 \ge  \sqrt{\frac{\beta(t_0)}{\alpha(t_0)}}, \phantom{a} \gamma_n \in (-\infty; + \infty), \phantom{a} n = \overline{1,3},$  Eq. (1.2) has a solution  $q(t)\equiv q_0(t) - i q_1(t) - j q_2(t) - k q_3(t)$ on $[t_0; + \infty)$ with $q_n(t_0) = \gamma_n, \phantom{a} n = \overline{0,3}$, and
$$
q_0(t) \ge  \sqrt{\frac{\beta(t)}{\alpha(t)}}, \phantom{aaa} t\ge t_0.
$$
\phantom{aaaaaaaaaaaaaaaaaaaaaaaaaaaaaaaaaaaaaaaaaaaaaaaaaaaaaaaaaaaaaaaaaaaaaaaaa} $\Box$}

   \vskip 10pt

{\bf 4. The case when $a_0(t)$ may change sign.} In this section we consider the case when $a_0(t)$ my change sign.
 Set:
$$
\left[\frac{\sqrt{\sum_{n=1}^3(b_n(t) + c_n(t))^2}}{ a_0(t)}\right]_0 \equiv \sist{\frac{\sqrt{\sum_{n=1}^3(b_n(t) + c_n(t))^2}}{ a_0(t)}, \ph if \ph a_0(t) \ne 0;}{ \phh  0, \phh \ph  if \ph a_0(t) = 0,} \phantom{aaaaaaaaaaaaaaaa}
$$
$$
\mathfrak{M}(t) \equiv \il{t_0}{t}||(d_1(\tau), d_2(\tau), d_3(\tau))|| d\tau + \frac{1}{2} \sup\limits_{\tau \in [t_0;t]}\left[\frac{\sqrt{\sum_{n=1}^3(b_n(\tau) + c_n(\tau))^2}}{ a_0(\tau)}\right]_0, \phantom{aaa}
$$
$$
\phantom{aaaaaaaaaa} R_\Gamma(t) \equiv |a_0(t)|(\Gamma + \mathfrak{M}(t))^2 + \sum\limits_{n=1}^3|b_n(t) + c_n(t)|(\Gamma + \mathfrak{M}(t)), \phh t\ge t_0,
$$
where $\Gamma> 0$ is a parameter. For any quaternion $q \equiv q_0 + i q_1 + j q_2 +k q_3 \ph (q_n \in \mathbb{R}, \ph n= \overline{0.3}),$ set $[q]_v \equiv (q_1, q_2, q_3)$.

{\bf Theorem 4.1.} {\it Let $\alpha_m(t)$ and $\beta_m(t), \ph m=1,2$ be the same as in Theorem 2.5. If:

\noindent
1) $a_n(t) \equiv 0, \ph n= \overline{1,3};$

\noindent
2) $\alpha_1(t) \le a_0(t) \le \alpha_2(t), \phh \beta_1(t) \le R_\Gamma(t) + d_0(t) \le \beta_2(t), \ph t\in [t_0;\tau_0);$

\noindent
3)   $b_0(t) + c_0(t) \ge  \frac{1}{2}\biggl(\frac{\alpha'_m(t)}{\alpha_m(t)} - \frac{\beta'_m(t)}{\beta_m(t)}\biggr) +
 2(-1)^m \sqrt{\alpha_m(t)\beta_m(t)}, \ph t\in [t_0;\tau_0), \ph m=1,2;$

\noindent
4) $b_0(t) + c_0(t) \ge 2|a_0(t)| R_\Gamma(t), \phh t\in [t_0;\tau_0);$

\noindent
5) $supp \hskip 3pt (b_n(t) + c_n(t)) \subset supp \hskip 3pt a_0(t), \ph n =\overline{1,3}$, the function $\left[\frac{\sqrt{\sum_{n=1}^3(b_n(t) + c_n(t))^2}}{ a_0(t)}\right]_0$  is bounded on $[t_0;\tau_0)$ if $\tau_0 < +\infty$ and is locally bounded on $[t_0;\tau_0)$ if $\tau_0 = + \infty$,

\noindent
then for every $\gamma_0 \in \biggl[- \sqrt{\frac{\beta_2(t_0)}{\alpha_2(t_0)}}; \sqrt{\frac{\beta_1(t_0)}{\alpha_1(t_0)}}\biggr], \ph \gamma_n \in \mathbb{R}, \ph n = \overline{1,3},$
with $||(\gamma_1, \gamma_2, \gamma_3)|| \le \Gamma$ Eq.~ (1.1) has a solution $q(t) \equiv q_0(t) - i q_1(t) - j q_2(t) - k q_3(t)$ on $[t_0;\tau_0)$ satisfying the initial conditions  $q_n(t_0) =~ \gamma_n, \ph  n =~ \overline{0,3}$, and
$$
- \sqrt{\frac{\beta_2(t)}{\alpha_2(t)}} \le q_0(t) \le \sqrt{\frac{\beta_1(t)}{\alpha_1(t)}}, \phh t\in [t_0;\tau_0); \eqno (4.1)
$$
$$
||[q(t)]_v|| \le ||[q(t_0)]_v|| + \mathfrak{M}(t), \phh t\in [t_0;\tau_0). \eqno (4.2)
$$

If $\tau_0 < +\infty$ then $q(t)$ is continuable on $[t_0;\tau_0].$
}

Proof.
 Let $q(t) \equiv q_0(t) - i q_1(t) - j q_2(t) - k q_3(t)$ be the solution of Eq. (1.1) with $q_n(t_0) = \gamma_n, \ph n = \overline{0,3},$ and let $[t_0;T)$ be the maximum existence interval for $q(t)$. We must show that
$$
T \ge \tau_0. \eqno (4.3)
$$
Under the restriction 1) the system (2.1) takes the form
$$
\sist{q_0' + a_0(t) q_0^2 + \{b_0(t) + c_0(t)\} q_0 - P(t, q_1,q_2,q_3) =0;}{\widetilde{q} \hskip 1.5pt' + \mathcal{L}_{q_0}(t)\widetilde{q} - f_{q_0}(t) = 0, \phh t\ge t_0,} \eqno (4.4)
$$
where $f_{q_0}(t) \equiv ((b_1(t) + c_1(t)) q_0 + d_1(t), (b_2(t) + c_2(t)) q_0 + d_2(t), (b_3(t) + c_3(t)) q_0 + d_3(t))$,
$$
\mathcal{L}_{q_0}(t) \equiv \begin{pmatrix} b_0(t) + c_0(t) + 2 a_0(t) q_0 & c_3(t) - b_3(t) & b_2(t)  -c_2(t)\\
b_3(t) - c_3(t) & b_0(t) + c_0(t) + 2 a_0(t) q_0 & c_1(t) - b_1(t) \\
c_2(t) - b_2(t) & b_1(t) - c_1(t) & b_0(t) + c_0(t) + 2 a_0(t) q_0
\end{pmatrix},
$$
$t\ge t_0,    \ph \widetilde{q} \equiv (q_1, q_2, q_3).$
Since the hermitian part $\mathcal{L}_{q_0(t)}^H(t)$ of the matrix  $\mathcal{L}_{q_0(t)}(t)$ is $\mathcal{L}_{q_0(t)}^H(t) = diag \{ b_0(t) + c_0(t) + 2 a_0(t) q_0(t), b_0(t) + c_0(t) + 2 a_0(t) q_0(t), b_0(t) + c_0(t) + 2 a_0(t) q_0(t)\},$ by the second equation of the system (4.4)  $||[q(t)]_v||$ we have the estimate (see  [Hartm, p. 56, Lemma 4.2]):
$$
||[q(t)]_v|| \le ||[q(t_0)]_v||\exp\biggl\{-\il{t_0}{t}(b_0(\tau) + c_0(\tau) + 2 a_0(\tau) q_0(\tau)) d\tau\biggr\}+ \phantom{aaaaaaaaaaaaaaaaaaaaaaaaaaaaaaaa}
$$
$$
+ \il{t_0}{t}\exp\biggl\{-\il{\tau}{t}(b_0(s) + c_0(s) +2 a_0(s) q_0(s)) d s\biggr\}||f_{q_0(\tau)}(\tau)|| d \tau, \ph t\in [t_0;t_1). \eqno (4.5)
$$
From the condition 4) of the theorem it follows that
$$
b_0(t) + c_0(t) + 2 a_0(t) q_0(t) \ge 0, \phh t \in [t_0;t_1), \eqno (4.6)
$$
for some $t_1 > t_0.$ Show that
$$
-\sqrt{\frac{\beta_2(t)}{\alpha_2(t)}} \le q_0(t) \le  \sqrt{\frac{\beta_1(t)}{\alpha_1(t)}}, \phh t\in [t_0; t_1); \eqno (4.7)
$$
$$
||[q(t)]_v|| \le ||[q(t_0)]_v|| + \mathfrak{M}(t), \phh t\in [t_0;t_1). \eqno (4.8)
$$
From (4.5) and (4.6) it follows
$$
||[q(t)]_v|| \le ||[q(t_0)]_v|| + \frac{1}{2}\exp\biggl\{-\il{t_0}{t}(b_0(s) + c_0(s) + 2 a_0(s) q_0(s)) ds \biggr\} \times \phantom{aaaaaaaaaaaaaaaaaaa}
$$
$$
\times \il{t_0}{t}\biggl(\exp\biggl\{\il{t_0}{\tau}(b_0(s) + c_0(s) + 2a_0(s) q_0(s))d s\biggr\}\biggr)'\left[\frac{\sqrt{\sum_{n=1}^3(b_n(\tau) + c_n(\tau))^2}}{ a_0(\tau)}\right]_0 d \tau + \phantom{aaaaa}
$$
$$
\phantom{aaaaaaaaaaaaaaaaaaaaaaaaaaaaaaaaa} + \il{t_0}{t}||(d_1(\tau), d_2(\tau), d_3(\tau))|| d\tau, \ph t\in [t_0;t_1).
$$
From here from (4.6) and 5) it follows (4.8). Since $||[q(t_0)]_v|| \le \Gamma$ from (4.8) we obtain
$$
- R_\Gamma(t) + q_0(t) \le P(t, q_1(t), q_2(t), q_3(t)) \le  R_\Gamma(t) + q_0(t), \phh t\in [t_0;t_1).
$$
From here and from 2) it follows
$$
\beta_1(t) \le P(t, q_1(t), q_2(t), q_3(t)) \le \beta_2(t), \phh t\in [t_0;t_1). \eqno (4.9)
$$
Consider the Riccati equation
$$
r' + a_0(t) r^2 + \{b_0(t) + c_0(t)\} r  - P(t, q_1(t), q_2(t), q_3(t)) = 0, \phh t\in [t_0;t_1). \eqno (4.10)
$$
Let $r(t)$ be a solution of this equation with $r(t_0) = q_0(t_0)$. Then by virtue of Theorem 2.1  from   1), 2) and (4.9) it follows that $r(t)$ exists on $[t_0;t_1)$ and
$$
-\sqrt{\frac{\beta_2(t)}{\alpha_2(t)}} \le r(t) \le  \sqrt{\frac{\beta_1(t)}{\alpha_1(t)}}, \phh t\in [t_0; t_1).
$$
Obviously $q_0(t)$ is a solution of Eq. (4.7) on $[t_0;t_1)$. Hence by the uniqueness theorem $q_0(t)$ coincides with $r(t)$ on $[t_0;t_1)$, and therefore (4.7)  is valid.
Let $T_1$ be the upper bound of all $t_1 \in [t_0;T)$ for which (4.7) - (4.9)  are satisfied. We assert that
$$
T_1 = T. \eqno (4.11)
$$
Indeed otherwise from (4.7)  it follows that
$$
q_0(t) \ge -\sqrt{\frac{\beta_2(T_1)}{\alpha_2(T_1)}}.
$$
From here and from 4)  it follows that $b_0(t) + c_0(t) + 2 a_0(t) q_0(t) \ge 0, \ph t\in [T_1;T_2)$ for some $T_2 > T_1$. Hence
$$
b_0(t) + c_0(t) + 2 a_0(t) q_0(t) \ge 0, \ph t\in [t_0;T_2). \eqno (4.12)
$$
Then repeating the arguments of the proof of (4.7) and (4.8) we conclude that
$$
-\sqrt{\frac{\beta_2(t)}{\alpha_2(t)}} \le q_0(t) \le  \sqrt{\frac{\beta_1(t)}{\alpha_1(t)}}, \phh t\in [t_0;T_2);
$$
$$
||[q(t)]_v|| \le ||[q(t_0)]_v|| + \mathfrak{M}(t), \phh t\in [t_0:T_2),
$$
which with (4.12)  contradicts the definition of $T_1$. The obtained contradiction proves (4.11). Thus
$$
-\sqrt{\frac{\beta_2(t)}{\alpha_2(t)}} \le q_0(t) \le  \sqrt{\frac{\beta_1(t)}{\alpha_1(t)}}, \phh t\in [t_0;T);
$$
$$
||[q(t)]_v|| \le ||[q(t_0)]_v|| + \mathfrak{M}(t), \phh t\in [t_0:T).
$$
By virtue of Lemma 2.1 from here it follows (4.3) and fulfillment of (4.1) and (4.2). If $\tau_0 < + \infty$ then by Lemma 2.1 from (4.1) and (4.2) it follows that $q(t)$ is continuable on $[t_0;\tau_0]$. The theorem is proved.

\vskip 20pt
Let $\tau_0 < + \infty$. Set:

$$
\mathfrak{M}^*(t) \equiv \il{t}{\tau_0}||(d_1(\tau), d_2(\tau), d_3(\tau))|| d\tau + \frac{1}{2} \sup\limits_{\tau \in [t;\tau_0]}\left[\frac{\sqrt{\sum_{n=1}^3(b_n(\tau) + c_n(\tau))^2}}{ a_0(\tau)}\right]_0, \phantom{aaaaaaaaaaaaaaaaaaaaaaaaaaaaaaaaaaaaaaaaaaaaaaaa}
$$
$$
\phantom{aaaaaaaaaaaaaaaa} R_\Gamma^*(t) \equiv |a_0(t)|(\Gamma + \mathfrak{M}^*(t))^2 + \sum\limits_{n=1}^3|b_n(t) + c_n(t)|(\Gamma + \mathfrak{M}^*(t)), \phh t\in  [t_0;\tau_0].
$$

{\bf Corollary 4.1.} {\it Let $\alpha_m(t)$ and $\beta_m(t), \ph m=1,2,$ be continuously differentiable on $[t_0;\tau_0]$ functions such that $(-1)^m\alpha_m(t) > 0, \ph (-1)^m\beta(t) > 0, \ph t\in [t_0;\tau_0], \ph m=1,2.$ If:

\noindent
1) $a_n(t) \equiv 0, \ph n = \overline{1,3};$

\noindent
1$^*$)  $\alpha_1(t) \le a_0(t) \le \alpha_2(t), $

\noindent
2$^*$)  $b_0(t) + c_0(t) \le -\frac{1}{2}\biggl(\frac{\alpha'_m(t)}{\alpha_m(t)} - \frac{\beta'_m(t)}{\beta_m(t)}\biggr) + 2(-1)^m\sqrt{\alpha_m(t)\beta_m(t)}, \ph t\in [t_0;\tau_0], \ph m=1,2;$

\noindent
3$^*$) $b_0(t) + c_0(t) \le -2|a_0(t)| R_\Gamma^*(t), \ph t\in [t_0;\tau_0];$

\noindent
4$^*$) $supp \hskip 3pt (b_n(t) + c_n(t)) \subset supp \hskip 3pt a_0(t), \ph n =\overline{1,3}$, the function $\left[\frac{\sqrt{\sum_{n=1}^3(b_n(t) + c_n(t))^2}}{ a_0(t)}\right]_0$  is bounded on $[t_0;\tau_0],$

\noindent
then for every $\gamma_0 \in \biggl[- \sqrt{\frac{\beta_1(\tau_0)}{\alpha_1(\tau_0)}}; \sqrt{\frac{\beta_2(\tau_0)}{\alpha_2(\tau_0)}}\biggr], \ph \gamma_n \in \mathbb{R}, \ph n = \overline{1,3},$
with $||(\gamma_1, \gamma_2, \gamma_3)|| \le \Gamma$ Eq.~ (1.1) has a solution $q(t) \equiv q_0(t) - i q_1(t) - j q_2(t) - k q_3(t)$ on $[t_0;\tau_0]$ satisfying the initial conditions  $q_n(\tau_0) =~ \gamma_n, \ph  n =~ \overline{0,3}$, and
$$
- \sqrt{\frac{\beta_1(t)}{\alpha_1(t)}} \le q_0(t) \le \sqrt{\frac{\beta_2(t)}{\alpha_2(t)}}, \phh t\in [t_0;\tau_0]; \eqno (4.13)
$$
$$
||[q(t)]_v|| \le ||[q(\tau_0)]_v|| + \mathfrak{M}^*(t), \phh t\in [t_0;\tau_0]. \eqno (4.14)
$$
}

Proof. Set: $\lambda_0 \equiv t_0 = \tau_0, \ph \widetilde{a}(t) \equiv - a(\lambda_0 - t), \ph  \widetilde{b}(t) \equiv - b(\lambda_0 - t), \ph  \widetilde{c}(t) \equiv - c(\lambda_0 -~ t), \linebreak  \widetilde{d}(t) \equiv - d(\lambda_0 - t), \ph  \widetilde{a}_0(t) \equiv - a_0(\lambda_0 - t), \ph  \widetilde{b}_n(t) \equiv - b_n(\lambda_0 - t), \ph \widetilde{c}_n(t) \equiv - c_n(\lambda_0 -~ t), \linebreak \widetilde{d}_n(t) \equiv~ - d_n(\lambda_0 - t),$
$$
\widetilde{\mathfrak{M}}(t) \equiv \il{t_0}{t}||(\widetilde{d}_1(\tau), \widetilde{d}_2(\tau), \widetilde{d}_3(\tau))|| d\tau + \frac{1}{2} \sup\limits_{\tau \in [t_0;t]}\left[\frac{\sqrt{\sum_{n=1}^3(\widetilde{b}_n(\tau) + \widetilde{c}_n(\tau))^2}}{ \widetilde{a}_0(\tau)}\right]_0, \phantom{aaa}
$$
$$
\phantom{aaaaaaaaaa} \widetilde{R}_\Gamma(t) \equiv |\widetilde{a}_0(t)|(\Gamma + \widetilde{\mathfrak{M}}(t))^2 + \sum\limits_{n=1}^3|\widetilde{b}_n(t) + \widetilde{c}_n(t)|(\Gamma + \widetilde{\mathfrak{M}}(t)), \phh t\in [t_0;\tau_0],
$$
where
$$
\left[\frac{\sqrt{\sum_{n=1}^3(\widetilde{b}_n(t) + \widetilde{c}_n(t))^2}}{\widetilde{a}_0(t)}\right]_0 \equiv \sist{\frac{\sqrt{\sum_{n=1}^3(\widetilde{b}_n(t) + \widetilde{c}_n(t))^2}}{\widetilde{a}_0(t)}, \ph if \ph \widetilde{a}_0(t) \ne 0;}{ \phh  0, \phh \ph  if \ph \widetilde{a}_0(t) = 0.}
$$
In Eq. (1.1) make the substitution
$$
q(t) = u(\lambda_0 - t), \phh t\in [t_0;\tau_0].
$$
we obtain
$$
u' + u \widetilde{a}(t) u + \widetilde{b}(t) u + u \widetilde{c}(t) + \widetilde{d}(t) = 0, \phh t\in [t_0;\tau_0]. \eqno (4.15)
$$
It is not difficult to verify that
$$
\widetilde{\mathfrak{M}}(\lambda_0 - t) = \mathfrak{M}^*(t), \ph \widetilde{R}_\Gamma(\lambda_0 - t) = R^*_\Gamma(t), \ph t\in [t_0;\tau_0].
$$
From here and from the conditions 1), 1$^*$) - 4$^*$) of the corollary we get:
$$
\widetilde{\alpha}_1(t) \le \widetilde{a}_0(t)\le \widetilde{\alpha}_2(t),\ph \widetilde{\beta}_1(t) \le \widetilde{R}_\Gamma(t) + \widetilde{d}_0(t) \le \widetilde{\beta}_2(t), \ph \widetilde{b}_0(t) + \widetilde{c}_0(t) \ge 2|\widetilde{a}_0(t)| \widetilde{R}_\Gamma(t),
$$
$$
\widetilde{b}_0(t) + \widetilde{c}_0(t) \ge \frac{1}{2}\biggr(\frac{\widetilde{\alpha}'_m(t)}{\widetilde{\alpha}_m(t)} - \frac{\widetilde{\beta}'_m(t)}{\widetilde{\beta}_m(t)}\biggr) + 2(-1)^m\sqrt{\widetilde{\alpha}_m(t)  \widetilde{\beta}_m(t)}, \phh t\in [t_0;\tau_0],
$$
where $\widetilde{\alpha}_m(t) \equiv - \alpha_{3-m}(\lambda_0 - t), \ph \widetilde{\beta}_m(t) \equiv - \beta_{3-m}(\lambda_0 - t), \ph m=1,2, \ph t\in [t_0; \tau_0], \linebreak supp \hskip 3pt (\widetilde{b}_n(t) + \widetilde{c}_n(t)) \subset supp \hskip 3pt \widetilde{a}_0(t), \ph n= \overline{1,3},$ the function
$
\left[\frac{\sqrt{\sum_{n=1}^3(\widetilde{b}_n(t) + \widetilde{c}_n(t))^2}}{\widetilde{a}_0(t)}\right]_0
$
is bounded on $[t_0; \tau_0]$. By Theorem 4.1 from here is seen that for every $\gamma_0 \in \biggl[-\sqrt{\frac{\widetilde{\beta}_2(t_0)}{\widetilde{\alpha}_2(t_0)}}; \sqrt{\frac{\widetilde{\beta}_1(t_0)}{\widetilde{\alpha}_1(t_0)}}\biggr], \linebreak \gamma_n \in~ \mathrm{R}, \ph n= \overline{1,3},$ with $||(\gamma_1, \gamma_2, \gamma_3)|| \le \Gamma$ Eq. (4.15) has a solution $u(t) \equiv u_0(t) - i u_1(t) - j u_2(t) - k u_3(t)$ on $[t_0;\tau_0]$ and
$$
-\frac{\widetilde{\beta}_2(t)}{\widetilde{\alpha}_2(t)} \le u_0(t) \le \frac{\widetilde{\beta}_1(t)}{\widetilde{\alpha}_1(t)}, \phh ||[u(t)]_v|| \le  ||[u(t_0)]_v|| + \widetilde{\mathfrak{M}}(t),  \phh t\in [t_0;\tau_0],
$$
From here it follows that Eq. (1.1) has a solution $q(t) \equiv q_0(t) - i q_1(t) - j q_2(t) - k q_3(t)$ on $[t_0;\tau_0]$, satisfying the initial conditions  $q_n(\tau_0) = \gamma_n, \ph n = \overline{0,3}$ and the estimates (4.13) and (4.14) are valid. The corollary is proved.

\vskip 10 pt

{\bf 5. A completely non conjugation theorem}. Consider the linear system
$$
\left\{
\begin{array}{l}
\phi' = C(t) \phi + A(t) \psi;\\
\phantom{aaa}\\
\psi' = - D(t) \phi - B(t) \psi, \phantom{aaa}t \ge t_0.
\end{array}
\right.
\eqno (5.1)
$$
where $\phi = \phi(t)$   and  $\psi = \psi(t)$  are the unknown continuously differentiable vector functions of dimension 4,   $A(t), \phantom{a} B(t), \phantom{a} C(t)$ and  $D(t)$  are the same matrix functions as in (2.5).

{\bf Definition 5.1.} {\it  We will say that the solution  $(\phi(t), \psi(t))$ of the system (5.1) satisfies the completely non conjugation condition if  $\phi(t)\ne \theta, \phantom{a} \psi(t) \ne \theta \phantom{a} t\ge t_0,$ where $\theta$  is the null vector of dimension 4.}

{\bf Theorem  5.1.} {\it  Let the conditions of Theorem 3.1   (of Theorem 3.2) are satisfied. Then the solution   $(\phi(t), \psi(t))$  of the system (5.1) with   $\psi(t_0) = (\gamma_0 E - \gamma_1 I - \gamma_2 J - \gamma_3 K) \phi (t_0) \ne \theta$,  where  $\gamma_n \ge 0, \phantom{a} n  \in \mathfrak{S}(\ne\emptyset), \phantom{a} \sum\limits_{n\in\mathfrak{S}}\gamma_n \ne 0,\ph \gamma_n \in (-\infty;+\infty), \ph n\in \mathfrak{O}$     (where   $\gamma_n > 0, \linebreak n =~ \overline{0,3}$) satisfies of the completely non conjugation condition.}

Proof. Let the conditions of Theorem 3.1   (of Theorem 3.2) be satisfied and let
$q(t) \equiv q_0(t) - i q_1(t) - j q_2(t) - k q_3(t)$ be the solutions of Eq. (1.2) with  $q_n(t_0) = \gamma_0, \phantom{a} n= \overline{0,3}$
By virtue of Theorem 3.1   (Theorem 3.2) $q(t)$ exists on  $[t_0;+\infty)$.  From the condition  $\sum\limits_{n\in\mathfrak{S}} \gamma_n >~ 0 \phantom{a} (\gamma_n > 0, \phantom{a} n = \overline{0,3})$  it follows that
$$
q(t) \ne 0, \phantom{aaa} t\ge t_0. \eqno (5.2)
$$
By (2.4) $Y_1(t)\equiv \widehat{q(t)}$ is a solution of Eq. (2.3) on  $[t_0;+\infty)$.  From  (5.2) it follows that
$$
\det Y_1(t) \ne 0, \phantom{aaa} t\ge t_0. \eqno (5.3)
$$
Let  $\Phi_1(t)$  be the solution of the matrix equation
$$
\Phi' = [A(t) Y_1(t) + C(t)] \Phi = 0, \phantom{aaa} t\ge t_0,
$$
satisfying the initial condition $\Phi_1(t_0)= E$. Them by the LIouville's formula we have
$$
\det \Phi(t) = \exp\biggl\{\int\limits_{t_0}^t tr[A(\tau) Y_1(\tau) + C(\tau)] d\tau\biggr\} > 0, \phantom{aaa} t\ge t_0. \eqno (5.4)
$$
Let $(\phi(t), \psi(t))$ be the solution of the system (5.1) satisfying the initial condition of the theorem. Then
$$
\phi(t) = \Phi(t) \phi(t_0), \phantom{aaa} \psi(t) = Y_1(t) \Phi(t) \phi(t_0).
$$
From here from (5.3) and  (5.4)  it follows that  $\phi(t) \ne \theta, \phantom{a} \psi(t) \ne \theta, \phantom{a} t\ge t_0$. The theorem is proved.

{\bf Remark  5.1.} {\it Except in a special case when  $A(t)$ and  $D(t)$   are diagonal matrices and  $C(t) = B^*(t), \phantom{a} t\ge t_0$ (here $*$  is the transpose sign) the system (5.1) is not hamiltonian.}

\vskip 20pt

 \centerline{\bf References}

\vskip 20pt

\noindent
1. J. D. Gibbon, D. D. Halm, R. M. Kerr and I. Roulstone. Quatrnions and particle \linebreak \phantom{aa} dynamics in the Euler fluid equations. Nonlinearity, vol. 19, pp. 1962 - 1983, 2006.

\noindent
2. V. Christiano and F. Smarandache. An Exact Mapping from Navier - Stokes Equation\linebreak \phantom{a} to Schrodinger Equation via Riccati Equation. Progress in Phisics, vol. 1,\linebreak \phantom{aa} pp. 38, 39, 2008.

\noindent
3. H. Zoladek. Classifications of diffeomormhisms of $S^4$ induced by quaternionic Riccati \linebreak \phantom{aa} equations with periodic coefficients. Topological methods in Nonlinear Analysis. Journal\linebreak \phantom{aa} of the Juliusz Schauder Center, vol. 33, pp. 205 - 215, 2009.

\noindent
4. K. Leschke and K. Morya. Application of Quaternionic Holomorphic Geometry to\linebreak \phantom{aa} minimal surfaces. Complex Manifolds, vol. 3, pp. 282 - 300, 2006.

\noindent
5. J. Campos, J. Mawhin. Periodic solutions of quaternionic - valued ordinary \linebreak \phantom{aa} differential equations, Anbnali di Mathematica, vol. 185, pp. 109 -  127, 2006.

\noindent
6. P Wilzinski, Quaternionic - valued differential equations. The Riccati equation. Journal\linebreak \phantom{aa} of Differential Equations, vol. 247, pp. 2163 - 2187, 2009.

\noindent
7. G. A. Grtigorian. Global solvability criteria for scalar Riccati equations with complex \linebreak \phantom{aa}  coefficients. Differ. Uravn. 53 (2017) no. 4m, pp. 459 - 464.

\noindent
8. Ph. Hartman, Ordinary differential equations, SIAM - Society for industrial and\linebreak \phantom{aaa} applied Mathematics, Classics in Applied Mathematics 38, Philadelphia 2002.

\noindent
9. G. A. Grigorian.  On two comparison tests for second-order linear  ordinary differential\linebreak \phantom{aa}  equations (Russian) Differ. Uravn. 47 (2011), no. 9, 1225 - 1240; translation in \linebreak \phantom{aa} Differ. Equ. 47 (2011), no. 9 1237 - 1252, 34C10.

\noindent
10. G. A. Grigorian, "Two Comparison Criteria for Scalar Riccati Equations with\linebreak \phantom{aa} Applications". Russian Mathematics (Iz. VUZ), 56, No. 11, 17 - 30 (2012).

\noindent
11. G. A. Grigorian, Global Solvability of Scalar Riccati Equations. Izv. Vissh.\linebreak \phantom{aa} Uchebn. Zaved. Mat., 2015, no. 3, pp. 35 - 48.

 \end{document}